\documentclass[reqno,11pt]{amsart}
\usepackage{bm,cancel}
\usepackage{amssymb,amsfonts}
\usepackage{amsmath,latexsym}
\usepackage{xcolor}
\usepackage{tikz}
\usetikzlibrary{arrows.meta}
\usepackage{float}
\usepackage{pdfsync}
\usepackage{mathrsfs}
\usepackage[numbers,sort&compress]{natbib}

\usepackage[unicode,breaklinks=true,colorlinks=true]{hyperref}
\hypersetup{
citecolor=[HTML]{006666},
linkcolor=[HTML]{8B0000},
urlcolor =[HTML]{800000}
}

\newtheorem{thm}{Theorem}[section]

\newtheorem{lemma}[thm]{Lemma}
\newtheorem{prop}[thm]{Proposition}
\theoremstyle{definition}

\newtheorem{rem}[thm]{Remark}
\numberwithin{equation}{section}

\newcommand{\mR}{\mathbb R}

\newcommand{\mD}{\m{D}}

\newcommand{\dx}{\,{\rm d}x}

\newcommand{\bn}{\bar{\nabla}}
\newcommand{\dr}{\,{\rm d}r}

\newcommand{\al}{\alpha}

\newcommand{\f}{\frac}
\newcommand{\q}{\quad}
\newcommand{\na}{\nabla}
\newcommand{\un}{\underbrace}
\newcommand{\vp}{\varphi}

\renewcommand{\div}{\mathop{\rm div}\nolimits}

\renewcommand{\leq}{\leqslant}
\renewcommand{\geq}{\geqslant}
\renewcommand{\d}{\,{\rm d}}

\def\XXint#1#2#3{{\setbox0=\hbox{$#1{#2#3}{\int}$} \vcenter{\hbox{$#2#3$}}\kern-.5\wd0}}

\newcommand{\p}{\partial}

\newcommand{\esssup}{{\rm ess}\sup}
\def\be{\begin{equation}}
\def\ee{\end{equation}}
\def\bes{\begin{equation*}}
\def\ees{\end{equation*}}
\def\bs{\begin{split}}
\def\es{\end{split}}
\def\bali{\begin{aligned}}
\def\eali{\end{aligned}}

\def\mR{\mathbb{R}}
\def\mS{\mathbb{S}}

\def\un{\underbrace}
\def\ol{\overline}
\def\al{\alpha}

\def\la{\lambda}

\def\Dl{\Delta}

\def\Si{\Sigma}
\def\i{\infty}
\def\p{\partial}
\def\f{\frac}
\def\na{\nabla}

\def\O{\Omega}
\def\s{\sqrt}
\def\q{\quad}
\def\qq{\qquad}

\def\m{\mathcal}

\allowdisplaybreaks[4]
\def\cd{\cdot}

\def\bm{\boldsymbol}
\def\dx{\,\mathrm{d}x}
\def\d{\,\mathrm{d}}
\def\ba{\begin{equation}\begin{aligned}}
\def\ea{\end{aligned}\end{equation}}
\def\bn{\[\begin{aligned}}
\def\en{\end{aligned}\]}
\def\l{\label}
\def\ed{\buildrel\hbox{\footnotesize def}\over =}
\def\ef{\eqref}

\begin{document}

\title[Flux bounds for Leray's problem with Navier slip]{Refined critical flux bounds for Leray's problem in a pipe with the Navier-slip boundary}
\author[Z. Li]{Zijin Li}
\address[Z. Li]
{School of Mathematics and Statistics, Nanjing University of Information Science and Technology, Nanjing
	210044, China}
\email{zijinli@nuist.edu.cn}

\author[Z. Wu]{Zhipeng Wu}
\address[Z. Wu]{School of Mathematics and Statistics,
Nanjing University of Information Science and Technology,
Nanjing 210044, China}
\email{wuzhipeng@nuist.edu.cn}
\date{\today}

\begin{abstract}
We study uniqueness for stationary incompressible Navier--Stokes flows
in infinite pipes with the Navier-slip boundary.
For a cylindrical‌ pipe of unit radius and a positive friction ratio, uniqueness
holds when $|\Phi|\leq2\pi\s{40/7}\approx15$. This significantly improves upon the result of the first author, Liu and Zhou (J. Differential Equations 461:114108,
2026) where the uniqueness requires $|\Phi|\leq\pi/16\approx 0.2$.
For small positive friction ratio $\alpha$, this flux condition can be improved to
$|\Phi|\leq4\s2\,\pi^2/\al$. For an infinite pipe with a bounded, smooth, simply
connected cross-section $\Sigma$, we prove uniqueness for every
$0<\al<\i$ provided that $|\Phi|\leq c_\Sigma(1+\al)/\al$,
where $c_\Sigma>0$ depends only on $\Sigma$.

\end{abstract}
\maketitle


\tableofcontents
\section{Introduction}

\subsection{The problem and the main results}
We consider the three-dimensional stationary Navier--Stokes equations
\be\l{eq:NS}
 \left\{\bali&
 -\Dl\bm{u}+(\bm{u}\cdot\na)\bm{u}+\na\pi=0\,,\\
 &\div\bm{u}=0\,,
 \eali\right.\q\text{in }\mD
\ee
with the Navier-slip boundary condition:
\be\l{eq:Navier}
 \bm{u}\cdot\bm n=0,\q
 2(\mS\bm{u}\cd\bm n)_{\tan}+\al\bm{u}_{\tan}=0
 \q\text{on }\p\mD.
\ee
Here $\bm{u}$ is the velocity, $\pi$ is the scalar pressure, and $\al\geq0$ is the friction ratio. The viscosity is normalized to $1$, and $\mD=\Sigma\times\mR$ is the infinite cylinder, where $\Sigma$ is a bounded, smooth, simply connected cross-section. We write $\mD_R=\Si\times(-R,R)$ for $R>0$. \(\bm n\) denotes the outward unit normal on the lateral boundary, and 
\bes
 \mS\bm{u}\ed\f12\bigl(\na\bm{u}+(\na\bm{u})^{\mathsf T}\bigr),
 \q \bm{u}_{\tan}\ed\bm{u}-(\bm{u}\cdot\bm n)\bm n.
\ees

The prescribed flux is
\be\l{eq:flux}
 \int_\Sigma u_3(x_h,x_3)\d x_h=\Phi.
\ee
It is independent of \(x_3\) by incompressibility and the normal
boundary condition in \eqref{eq:Navier}. 

A divergence-free field $\bm u\in H^1(\Si\times(-R,R))^3$ for
every finite $R>0$ is a weak solution of
\eqref{eq:NS}--\eqref{eq:Navier}--\eqref{eq:flux} if it satisfies
the normal boundary and flux conditions in the trace sense and
\bes\l{DEweak}
 2\int_{\mD}\mS\bm{u}:\mS\bm\vp\dx
 +\int_{\mD}((\bm{u}\cdot\na)\bm{u})\cdot\bm\vp\dx
 +\al\int_{\p\mD}\bm{u}_{\tan}\cdot\bm\vp_{\tan}\d S=0
\ees
for every smooth divergence-free test field $\bm\vp$ with compact
axial support and zero normal trace on the lateral boundary.
The colon denotes the Frobenius inner product. The pressure is
determined up to an additive constant.

We also include the no-slip case, denoted by \(\al=\i\).
In this case, \eqref{eq:Navier} is replaced by the homogeneous
Dirichlet condition
\[
\bm u=0\qquad\text{on }\p\mD.
\]
The weak formulation above is understood without the boundary integral, with the test fields required to have zero trace
on the lateral boundary.

If $\Sigma=B\ed\{x_h\in \mathbb{R}^2\,:\,|x_h|<1\}$, problem \eqref{eq:NS}--\eqref{eq:Navier}--\eqref{eq:flux} admits the generalized Hagen--Poiseuille flow:
\be\l{eq:Poiseuille}
 \bali
 \bm P_{\Phi,\text{disk}}&=P_{\Phi,\text{disk}}\bm e_3,\qq
 P_{\Phi,\text{disk}}(x_h)=\f{2\Phi}{(\al+4)\pi}
                  \bigl(\al+2-\al|x_h|^2\bigr),\q\pi_P(x)=-\f{8\al\Phi}{(\al+4)\pi}x_3.
 \eali
\ee
See~\cite[(1.11)]{LiLiuZhou}. For \(\al=\i\), the formulas in \eqref{eq:Poiseuille}
are understood as their limits:
\[
P_{\Phi,\text{disk}}(x_h)
=\f{2\Phi}{\pi}(1-|x_h|^2),
\qquad
\pi_P(x)=-\f{8\Phi}{\pi}x_3.
\]
Thus a fundamental question naturally arises:

\[
\boxed{\text{Is \eqref{eq:Poiseuille} the only solution of \eqref{eq:NS}--\eqref{eq:Navier}--\eqref{eq:flux}?}}
\]

Li, Liu and Zhou \cite[Section 4.4]{LiLiuZhou} proved uniqueness
for $\Phi\leq\pi/16\approx0.2$. For unit radius and viscosity,
this corresponds to a Reynolds number of order $10^{-1}$.
Experimental pipe-flow transition occurs at Reynolds numbers of
order $10^3$; see \cite{Binnie,Avila}. It is well known that the onset of turbulence can give rise to flow states beyond the laminar Poiseuille profile. The considerable difference between the Reynolds-number range covered by our uniqueness result and that observed in experiments highlights a gap of several orders of magnitude between current mathematical theory and experimental observations. Motivated by this gap, our first result improves the existing theoretical estimate by increasing the explicit uniform flux bound from $\pi/16$ to a value exceeding $15$.

\begin{thm}[A uniform flux bound]\l{thm:main}
Let $\Si=B$ and \(0<\al\leq\i\). Suppose that $\bm u$ is a weak solution of
\eqref{eq:NS}-\eqref{eq:Navier}-\eqref{eq:flux} in the sense defined above and satisfies
\be\l{eq:growth}
 \|\na\bm{u}\|_{L^2(\mD_\zeta)}=o(\zeta^{3/2})\,,
 \q\text{as }\q\zeta\to\i.
\ee
If
\bes
 |\Phi|\leq2\pi\s{\f{40}{7}}\approx 15,
\ees
then \(\bm{u}=\bm P_{\Phi,\text{disk}}\).
\end{thm}

The sufficient flux bound in Theorem~\ref{thm:main}
\bes
2\pi\s{\f{40}{7}}
 \approx15.0197>15
\ees
holds for $0<\al\leq\i$. However, for small positive friction ratio, the sufficient flux threshold tends to infinity. 
This result holds even for the cross section just being a simply connected smooth domain $\Sigma$. In this situation, the corresponding Poiseuille flow is the solution of the following boundary value problem:

\[
\left\{
\begin{aligned}
	&\bm P_\Phi = P_\Phi\bm e_3,\\[6pt]
	&-\Dl_h P_\Phi(x_h) = g_{\al,\Phi},
	 \quad\quad&&\text{in }\Si,\\[6pt]
	&\p_{\bm n_h}P_\Phi = -\al P_\Phi,
	 &&\text{on }\p\Si,\\[6pt]
	&\int_\Si P_\Phi(x_h)\d x_h = \Phi.
\end{aligned}
\right.
\]
Here $g_{\al,\Phi}$ is a constant fixed by the flux, $\bm n_h$ is the outward unit normal to $\p\Si$, and $\pi_P=-g_{\al,\Phi}x_3$ is the associated pressure. The dependence of $P_\Phi$ on $\al$ is left implicit. See \cite[(1.6)]{LiLiuZhou}.

\begin{thm}[Small positive friction]\l{thm:small-friction}
Let $\Si\subset\mR^2$ be bounded, smooth and simply connected,
and let {$0<\al<\i$}. Suppose that $\bm u$ is a weak solution
of \ef{eq:NS}--\ef{eq:flux} in $\mD=\Si\times\mR$ and satisfies
\ef{eq:growth}. There is a constant $c_\Si>0$, depending only on
$\Si$, such that
\bes
 |\Phi|\leq c_\Si\f{1+\al}{\al}
\ees
implies $\bm u=\bm P_\Phi$.
In particular, if $\Si=B$, then $\bm u=\bm P_{\Phi,\text{disk}}$
whenever
\be\l{eq:small-friction-threshold}
 |\Phi|<\f{\s2\,\pi^2(\al+4)}{\al}.
\ee
\end{thm}

\begin{rem}
\sl{With unit pipe radius and viscosity, the \emph{Reynolds number} in the pipe $\mD$ is
\[
\mathrm{Re}_{\mD}=\f{2|\Phi|}{\pi}.
\]
Thus by Theorem \ref{thm:small-friction},
\[
\mathrm{Re}_{\mD}\leq\f{8\s2\,\pi}{\al}
\]
guarantees uniqueness. The upper endpoint tends to infinity as \(\al\to0_+\). To see the role of friction, write the axial profile as
\bes
 P_\Phi(r)=\f\Phi\pi
       +\f{\al\Phi}{(\al+4)\pi}(1-2r^2).
\ees
At \(\Phi=4\s2\,\pi^2/\al\),
\bes
 \f\Phi\pi=\f{4\s2\,\pi}{\al}\rightarrow\i,
 \qq
 \p_rP_\Phi(r)=-\f{16\s2\,\pi}{\al+4}r
              \rightarrow-4\s2\,\pi r.
\ees
Thus the mean speed grows without bound, while the shear stays bounded.
The constant mean \(\Phi/\pi\) disappears upon differentiation and
does not enter the bulk interaction term in the energy estimate. For example, at \(\al=0.01\) and this endpoint flux,
\bes
 \mathrm{Re}_{\mD}\approx3554.3064,
 \qq
 \f{P_\Phi(1)}{\Phi/\pi}=\f4{4.01}\approx0.9975.
\ees
The wall speed is about \(99.75\%\) of the mean speed, whereas the Dirichlet condition imposes zero velocity at the wall.}
\end{rem}

\subsection{Related works}
Leray's problem concerns stationary viscous flows in an unbounded channel
or pipe with a prescribed flux. A basic question is whether the flux and
the boundary condition determine the velocity uniquely. For the no-slip
condition, the existence and uniqueness theory in unbounded domains was
developed by Amick~\cite{Amick} and Ladyzhenskaya and
Solonnikov~\cite{LadyzhenskayaSolonnikov}. See also
Galdi~\cite{Galdi} for a systematic introduction. The behavior at the ends of
the pipe is an important part of the problem, since a local weak formulation alone does
not prescribe the growth of the solution at infinity.

The Navier-slip boundary condition relates the tangential viscous stress tensor to the
tangential velocity on the boundary. Its origin goes back to Navier~\cite{Navier}.
For the stationary Stokes and Navier--Stokes equations with this boundary
condition, see Acevedo Tapia, Amrouche, Conca and Ghosh~\cite{Acevedo}.
In an infinite pipe, the friction coefficient changes not only the
Hagen--Poiseuille profile but also the boundary contribution in the energy
identity. An estimate uniform in that coefficient must therefore treat
these two effects together.

Li, Pan and Yang~\cite{LiPanYang} studied existence, uniqueness,
regularity and the behavior at infinity for small-flux solutions in
distorted pipes with Navier slip. In that result the smallness threshold
can depend on the friction coefficient. Wang and Xie~\cite{WangXie}
obtained uniqueness and uniform structural stability results for
axisymmetric Poiseuille flows. The present paper concerns general
three-dimensional weak solutions, without an axisymmetry assumption,
and explicit flux bounds with quantified dependence on the friction
coefficient.

Recently, Li, Liu and Zhou~\cite{LiLiuZhou} proved a refined uniqueness theorem for
Leray's problem with positive Navier friction. They assumed that the
\(L^2\)-norm of the velocity gradient on a pipe segment of length
\(2\zeta\) is \(o(\zeta^{3/2})\). Under this growth condition, they
obtained a positive small-flux uniqueness threshold that can be chosen
independently of the friction coefficient. For the unit disk, their
explicit uniform lower bound is \(\pi/16\approx0.2\). Their proof combines
cross-sectional inequalities, a pressure estimate on unit end slabs,
and a Saint-Venant argument.

In this paper we improve the explicit uniform bound for the unit disk
to \(2\pi\s{40/7}>15\). We use the same growth assumption and the same
end-slab pressure estimate. The main improvement is in the way we handle the bulk energy estimate.
First, subtracting the cross-sectional mean from the Poiseuille profile
gives a smaller coefficient in the interaction term. Second, when the
friction is larger than one, we retain the positive circumferential
boundary dissipation. See Section \ref{Section2.1} and Section \ref{Section2.2} below. By constructing a suitable positive cubic polynomial as a comparison function, we derive an explicit lower bound for the first Robin eigenvalue, thereby obtaining a sharper coercivity estimate that fully exploits the positive boundary dissipation. See Lemma \ref{lem:square} in Section \ref{Sec2}.

We also study the different question of how large the admissible flux
can be when the positive friction ratio $\al$ tends to zero. We prove
that, on the unit disk, the closed flux range \(|\Phi|\leq4\s2\,\pi^2/\al\) is
admissible for \(0<\al\leq1\). Thus the lower bound for the
fixed-friction uniqueness threshold tends to infinity. The proof
keeps the exact horizontal strain dissipation and uses the fact that
the Poiseuille shear is radial. A rigid rotation has zero horizontal
strain, but it also has zero radial component. This allows us to
estimate the interaction without controlling the whole horizontal
velocity by the strain first. Furthermore, for pipes with general smooth, bounded, and simply connected cross-sections, we establish a uniqueness criterion with an admissible flux of order $1/\alpha$ as $\alpha\to0_+$, showing that the divergence of the uniqueness threshold in the small-friction limit is not restricted to circular pipes.

Finally, we mention that for circular cross-sections, positive friction is essential for
uniqueness in this class. At zero friction, rigid rotations can be
superimposed on a constant axial flow without changing the flux;
see \cite[Remark 1.2]{LiLiuZhou} and \cite{LiPanYangTotalSlip}.

\subsection{Notations}

We use the following notation.

\begin{itemize}
	
	\item
	$C$ denotes a positive constant whose value may change from line to line, and $C_{a,b,\cdots}$ denotes a positive constant depending only on $a,b,\cdots$. We write $A\lesssim B$ if $A\leq CB$ for some $C>0$.
	
	\item
In the standard Euclidean basis $\bm e_1,\bm e_2,\bm e_3$,
we write $x=(x_h,x_3)$, where $x_h=(x_1,x_2)$,
and decompose a three-dimensional vector $\bm w$ as
$\bm w=w_1\bm e_1+w_2\bm e_2+w_3\bm e_3$,
with $\bm w_h=(w_1,w_2)$ denoting its horizontal component.
We use the notation
\bes
\na_h=(\p_1,\p_2),\q
\Dl_h=\p_1^2+\p_2^2,\q
\div_h\bm w_h=\p_1w_1+\p_2w_2.
\ees

	\item
	For $\bm\xi=(\xi_1,\xi_2)\in\mathbb R^2$, we write
	\[
	\bm\xi^\perp\ed(-\xi_2,\xi_1),
	\]
	which denotes the counterclockwise rotation of $\bm\xi$
	through an angle of $\pi/2$.

	\item
	$\mathfrak{h}=(h_1,h_2,h_3)$ denotes a multi-index, where
	$h_1,h_2,h_3\in\mathbb N\cup\{0\}$, and we write
	$|\mathfrak{h}|=h_1+h_2+h_3$ and
	$\nabla^{\mathfrak{h}}
	=\partial_{x_1}^{h_1}\partial_{x_2}^{h_2}\partial_{x_3}^{h_3}$.

	\item
	We use the standard Lebesgue and Sobolev spaces.
	For a domain \(\O\), the \(L^p\)-norm is
	\bes
	\|f\|_{L^p(\O)}
	\ed
	\begin{cases}
		\displaystyle
		\left(\int_{\O}|f(x)|^p\dx\right)^{1/p},
		&1\leq p<\i,\\[2mm]
		\displaystyle
		\esssup_{x\in\O}|f(x)|,
		&p=\i.
	\end{cases}
	\ees
	For vector- and matrix-valued functions, $|\cdot|$ denotes
	the Euclidean and Frobenius norms, respectively.
	\item
	$H^m$ denotes the $L^2$-based Sobolev space with its norm
	\[
	\|f\|_{H^m}
	\ed
	\sum_{0\leq|\mathfrak{h}|\leq m}
	\|\nabla^{\mathfrak{h}}f\|_{L^2}.
	\]
\end{itemize}

\subsection{Organization of the paper}
This paper is organized as follows: Section \ref{Sec2} and Section \ref{SecSec3} are devoted to the proofs of Theorem \ref{thm:main} and Theorem  \ref{thm:small-friction}, respectively.

\section{Proof of Theorem \ref{thm:main}}\l{Sec2}

\subsection{A scalar Robin inequality from a cubic polynomial}\l{Section2.1}

We first give an identity for a positive comparison function.
It applies to arbitrary scalar functions, without imposing a Robin
boundary condition on those functions.

\begin{lemma}\l{lem:square}
Let \(q\in C^2(\ol B)\) be positive, and suppose that
\be\l{BoundaryQ}
\p_{\bm n_h}q+(\al-1)q=0\q\text{ on}\q \p B,
\ee
 where
\(\bm n_h\) is the outward normal to the disk. For every
\(g\in H^1(B)\),
\be\l{eq:square}
 \bali
 &\int_B|\na_hg|^2\d x_h
       +(\al-1)\int_{\p B}|g|^2\d S_h=\int_B|\na_hg-g\na_h\log q|^2\d x_h
       +\int_B\f{-\Dl_hq}{q}|g|^2\d x_h.
 \eali
\ee
\end{lemma}

\begin{proof}
For smooth \(g\), the chain rule gives \(\na_h\log q=\na_hq/q\),
and expansion of the square gives
\ba\l{EI01001}
 \int_B|\na_hg-g\na_h\log q|^2\d x_h
 ={}&\int_B|\na_hg|^2\d x_h
 \un{-\int_B\f{\na_hq}{q}\cdot\na_h(g^2)\d x_h}_{I}+\int_B\f{|\na_hq|^2}{q^2}g^2\d x_h.
\ea
Integration by parts in the middle term yields
\ba\l{EII1001}
I=&-\int_{\p B}\f{\p_{\bm n_h}q}{q}g^2\d S_h+\int_B\left(\f{\Dl_hq}{q}
                 -\f{|\na_hq|^2}{q^2}\right)g^2\d x_h\\
=&(\al-1)\int_{\p B}g^2\d S_h+\int_B\left(\f{\Dl_hq}{q}
                 -\f{|\na_hq|^2}{q^2}\right)g^2\d x_h\,.              
\ea
Here in the second line we have applied the boundary condition of $q$ in \eqref{BoundaryQ}. Substituting \eqref{EII1001} in \eqref{EI01001}, the two terms containing \(|\na_hq|^2/q^2\) cancel, and we prove \ef{eq:square}. Finally, density and the trace theorem can extend the identity \eqref{eq:square} to \(H^1(B)\).
\end{proof}

In particular, if 
\[
-\Dl_hq\geq\la q,
\]
then \eqref{eq:square} implies
\be\l{eq:scalar-bound}
 \int_B|\na_hg|^2\d x_h
 + (\al-1)\int_{\p B}|g|^2\d S_h
 \geq\la\int_B|g|^2\d x_h
 \q\text{for all }g\in H^1(B).
\ee
We next construct \(q\) and \(\la\) explicitly for \(\al>1\). Here is the construction:

\begin{lemma}
{For every $\alpha >1, $} there exist a polynomial $q=q(x_h)$ and a constant $\la>0$ such that
\be\l{qcond}
\left\{
\begin{aligned}
&-\Dl_hq\geq\la q, \q&&\text{in}\q B\,;\\
&q>0, \q&&\text{in}\q B\,;\\
&\p_{\bm n_h}q+(\al-1)q=0,\q &&\text{on}\q \p B\,.\\
\end{aligned}
\right.
\ee
\end{lemma}

\begin{proof}
We seek $q$ in the following form
\be\l{Dpq}
 q=1+at+bt^2+ct^3,
\ee
where $t\ed 1-|x_h|^2$, and $a$, $b$, $c\in\mR$ are coefficients to be determined. 

We first focus on \eqref{qcond}$_3$. Since \(t=0\) and \(\p_{\bm n_h}t=-2\) on \(\p B\), we infer from \eqref{Dpq} that 
\bes
 q=1,\q \p_{\bm n_h}q=-2a\,.
\ees
So \eqref{qcond}$_3$ is satisfied by choosing
\bes
 a=\f{\al-1}{2}.
\ees

Noticing that \(\na_ht=-2x_h\) and \(\Dl_ht=-4\), direct calculation gives
\[
 -\Dl_ht=4,\q
 -\Dl_ht^2=-8+16t,\q
 -\Dl_ht^3=-24t+36t^2.
\]
Thus we have
\[
 (-\Dl_h-\la)q
 =(4a-8b-\la)+(16b-24c-\la a)t
       +(36c-\la b)t^2-\la ct^3.
\]
Choosing
\ba\l{Dbc}
 b=\f{4a-\la}{8},\q c=\f{16b-\la a}{24},
\ea
it follows that
\bes
 4a-8b-\la=0,\q16b-24c-\la a=0.
\ees
Consequently,
\be\l{eq:matching}
(-\Dl_h-\la)q=t^2\bigl((36-\la)c-\la b+\la c|x_h|^2\bigr).
\ee

To simplify the calculations, we make the following change of variables:
\be\l{eq:lambda-family}
 \mu\ed\f{\la(\al+2)}{\al-1},
\ee
and substituting this into \ef{Dbc}
gives
\be\l{eq:coefficients-K}
 \bali
 b&=\f{(\al-1)[2(\al+2)-\mu]}{8(\al+2)},\q c=\f{(\al-1)[(8-\mu)(\al-1)+24-4\mu]}{48(\al+2)}.
 \eali
\ee
Recalling \eqref{eq:matching}, to make $(-\Dl_h-\la)q\geq0$ in $B$, we need
\[
(36-\la)c-\la b\geq0.
\]
Substituting \eqref{Dbc} and \eqref{eq:lambda-family} in it, direct calculation shows
\ba\l{eq:K-condition}
 &(36-\la)c-\la b\\
 =&\f{\al-1}{48(\al+2)^2}\Big\{
 [36(\al+2)-\mu(\al-1)]\cdot[(8-\mu)(\al-1)+24-4\mu]\\
 &\hskip 3cm-6\mu(\al-1)[2(\al+2)-\mu]\Big\}\\
 =&\f{\al-1}{48(\al+2)^2}\bigl[(\mu^2-56\mu+288)(\al-1)^2+(10\mu^2-312\mu+1728)(\al-1)+2592-432\mu\bigr]\,.
\ea
Therefore, the coefficient of \((\al-1)^2\) in
\ef{eq:K-condition} is nonnegative if
\be\l{Krequest}
0<\mu\leq28-4\s{31}\approx 5.7289.
\ee
In such interval $\mu\in(0,28-4\s{31})$, it is clear that
\bes
 10\mu^2-312\mu+1728>0,\q2592-432\mu>0.
\ees
Thus it follows from
\ef{eq:matching}--\ef{eq:K-condition} that every such \(\mu\in(0,28-4\s{31})\)
gives a positive polynomial with \((-\Dl_h-\la) q\geq0\). This concludes \eqref{qcond}$_1$.

Finally, for \(\mu\in(0,28-4\s{31})\), both $b$ and $c$ in \eqref{eq:coefficients-K}
are positive when \(\al>1\), and thus \(q\geq1\) in \(\ol B\). This finishes the proof of \eqref{qcond}$_2$.

\end{proof}

By \eqref{Krequest}, we choose \(\mu=\f{40}{7}\) for computational convenience. In view of \eqref{eq:coefficients-K}
\bn
 b&=\f18\left(2(\al-1)
                   -\f{40(\al-1)}{7(\al+2)}\right)
   =\f{(\al-1)(7\al-6)}{28(\al+2)},\\
 c&=\f1{24}\left(16b
                   -\f{40(\al-1)}{7(\al+2)}\f{\al-1}{2}\right)
   =\f{(\al-1)(2\al-1)}{42(\al+2)}.
\en
Then the polynomial given in \eqref{Dpq} is
\[
 q_\al(x_h)={}1+\f{\al-1}{2}(1-|x_h|^2)+\f{(\al-1)(7\al-6)}{28(\al+2)}(1-|x_h|^2)^2+\f{(\al-1)(2\al-1)}{42(\al+2)}(1-|x_h|^2)^3.
\]

Recalling \eqref{eq:lambda-family}, one has the following result in view of \eqref{eq:scalar-bound}.

\begin{prop}\l{prop:Robin}
For every \(\al>1\) and every scalar \(g\in H^1(B)\),
\bes
 \int_B|\na_hg|^2\d x_h
 +(\al-1)\int_{\p B}|g|^2\d S_h
 \geq\f{40(\al-1)}{7(\al+2)}\int_B|g|^2\d x_h.
\ees
\end{prop}

\subsection{The perturbation and the bulk energy estimate}\l{Section2.2}
We first consider \(0<\al<\i\).
The Dirichlet case will be treated at the end of the
absorption argument.
Set \(\bm w=\bm{u}-\bm P_\Phi\) and
\(\bm w_h=(w_1,w_2)\). Since \(\bm{u}\) and \(\bm P_\Phi\)
have the same flux and satisfy the same boundary condition,
\be\l{eq:perturbation}
 \left\{
 \begin{aligned}
&-\Dl\bm w+(\bm w\cdot\na)\bm w
 +(\bm P_\Phi\cdot\na)\bm w
 +(\bm w\cdot\na)\bm P_\Phi+\na(\pi-\pi_P)&=0,\\
&\div\bm w=0,\qq
 \int_B w_3(x_h,x_3)\d x_h=0.&
\end{aligned}
\right.
\ee
On the lateral boundary, let \(\bm\tau\) be the unit
circumferential tangent, and write
\(w_\tau=\bm w\cdot\bm\tau\) and \(w_n=\bm w\cdot\bm n\).
The curvature of the unit circle is one, so \ef{eq:Navier} becomes
\be\l{eq:component-bc}
 w_n=0,\q \p_{\bm n}w_\tau=(1-\al)w_\tau,\q
 \p_{\bm n}w_3=-\al w_3.
\ee
Here \(w_\tau\) is a scalar component. In particular,
\(|\bm w_h|^2=|w_\tau|^2\) on the boundary.
The coordinate form \ef{eq:component-bc} is also the specialization
of~\cite[(1.5)]{LiLiuZhou} to the unit disk.

\subsubsection{Two forms of the viscous term}

Multiplying \ef{eq:perturbation} by \(\bm w\) on
\(\mD_\zeta\) gives the energy identity
\be\l{eq:energy}
 \bali
 \int_{\mD_\zeta}\bm w\cdot\Dl\bm w\dx
 ={}&\int_{\mD_\zeta}(\bm w_h\cdot\na_hP_\Phi)w_3\dx\\
 &+\left[\int_B\left\{\f12(w_3+P_\Phi)|\bm w|^2
                  +(\pi-\pi_P)w_3\right\}\d x_h
                                                    \right]_{-\zeta}^{\zeta}.
 \eali
\ee
The bracket denotes the value at \(x_3=\zeta\) minus the value
at \(x_3=-\zeta\). Indeed, incompressibility converts each
transport term into a boundary integral. The lateral contributions
vanish because \(\bm w\cdot\bm n=\bm P_\Phi\cdot\bm n=0\).
The pressure term gives the last term in the bracket. As in
\cite[Section 4.2]{LiLiuZhou}, these calculations can be obtained by
axial cutoffs in the weak formulation and then by integration over
the cut positions.

Using \ef{eq:component-bc}, integration by parts gives
\be\l{eq:gradient}
 \bali
 \int_{\mD_\zeta}\bm w\cdot\Dl\bm w\dx
 ={}&-\int_{\mD_\zeta}|\na\bm w|^2\dx
 -(\al-1)\int_{\p B\times(-\zeta,\zeta)}|w_\tau|^2\d S\\
 &-\al\int_{\p B\times(-\zeta,\zeta)}|w_3|^2\d S
 +\left[\int_B\bm w\cdot\p_3\bm w\d x_h
                                                    \right]_{-\zeta}^{\zeta}.
 \eali
\ee
There is also the symmetric-gradient identity
\be\l{eq:stress}
 \bali
 \int_{\mD_\zeta}\bm w\cdot\Dl\bm w\dx
 ={}&-2\int_{\mD_\zeta}|\mS\bm w|^2\dx
 -\al\int_{\p B\times(-\zeta,\zeta)}
                       (|w_\tau|^2+|w_3|^2)\d S\\
 &+\left[\int_B\sum_{i=1}^3
         w_i(\p_3w_i+\p_iw_3)\d x_h\right]_{-\zeta}^{\zeta}.
 \eali
\ee
To see this, use \(\Dl\bm w=2\div(\mathbb{S}\bm w)\), which follows
from \(\div\bm w=0\), and then apply \ef{eq:Navier} to
\(\bm w\). These are the two viscous identities used in
\cite[(4.5)--(4.6)]{LiLiuZhou}, with the curvature term kept
with its sign.

For \(0<\al<1\), calculating
\[
\al\times\ef{eq:gradient}+(1-\al)\times\ef{eq:stress},
\] 
the coefficients of the
circumferential boundary integral cancel, and the remaining bulk and lateral contributions are
\bes
 -\al\int_{\mD_\zeta}|\na\bm w|^2\dx
 -2(1-\al)\int_{\mD_\zeta}|\mS\bm w|^2\dx
 -\al\int_{\p B\times(-\zeta,\zeta)}|w_3|^2\d S.
\ees
For \(\al\geq1\), both lateral terms in \ef{eq:gradient}
are nonpositive. Consequently, for every \(\al>0\),
\be\l{eq:basic-dissipation}
 \bali
 \int_{\mD_\zeta}\bm w\cdot\Dl\bm w\dx
 \leq{}&-\min\{\al,1\}\int_{\mD_\zeta}|\na\bm w|^2\dx+C\int_{B\times\{\pm\zeta\}}
                         |\bm w|\,|\na\bm w|\d x_h.
 \eali
\ee
An integral over \(B\times\{\pm\zeta\}\) means the sum
over the two end sections. For larger \(\al\), we will use
the full identity \ef{eq:gradient} instead of discarding the
positive boundary dissipation.

\subsubsection{Centering the Poiseuille profile}

The profile \ef{eq:Poiseuille} satisfies
\be\l{eq:centered}
 P_\Phi-\f\Phi\pi
 =\f{\al\Phi}{(\al+4)\pi}(1-2|x_h|^2),
\ee
which follows that
\be\l{eq:centered1}
\left\|P_\Phi-\f\Phi\pi\right\|_{L^\i(B)}
 =\f{\al\Phi}{(\al+4)\pi}.
 \ee
The constant \(\Phi/\pi\) is both the mean and the midpoint of the range of \(P_\Phi\). Moreover, it corresponds to the Poiseuille flow when $\al=0$. Since
\(\div_h\bm w_h=-\p_3w_3\) and
\(\bm w_h\cdot\bm n_h=0\), horizontal integration by parts
gives
\be\l{eq:centered-interaction}
 \bali
 &\int_{\mD_\zeta}(\bm w_h\cdot\na_hP_\Phi)w_3\dx\\
 &=-\int_{\mD_\zeta}\left(P_\Phi-\f\Phi\pi\right)
                      \div_h(\bm w_hw_3)\dx\\
 &=-\int_{\mD_\zeta}\left(P_\Phi-\f\Phi\pi\right)
                      \bm w_h\cdot\na_hw_3\dx
 +\int_{\mD_\zeta}\left(P_\Phi-\f\Phi\pi\right)
                      w_3\p_3w_3\dx\\
 &=-\int_{\mD_\zeta}\left(P_\Phi-\f\Phi\pi\right)
                      \bm w_h\cdot\na_hw_3\dx
 +\f12\left[\int_B\left(P_\Phi-\f\Phi\pi\right)
                             |w_3|^2\d x_h\right]_{-\zeta}^{\zeta}.
 \eali
\ee
The additional end-section term is retained. By \ef{eq:centered},
its absolute value is at most
\be\l{eq:extra-end}
 \f{\al\Phi}{2(\al+4)\pi}
       \int_{B\times\{\pm\zeta\}}|w_3|^2\d x_h.
\ee

\subsubsection{Two estimates for the bulk interaction}

We first record a rough cross-sectional estimate valid for every
\(\al>0\). At almost every fixed \(x_3\), the divergence theorem
and the zero normal trace give
\bn
 \int_B|\bm w_h|^2\d x_h
 ={}&-\int_B(x_h\cdot\bm w_h)\div_h\bm w_h\d x_h
 -\sum_{i,j=1}^2\int_B w_i x_j\p_iw_j\d x_h\\
 \leq{}&\int_B|\bm w_h|
                \bigl(|\p_3w_3|+|\na_h\bm w_h|\bigr)\d x_h.
\en
Cauchy--Schwarz and \((a+b)^2\leq2a^2+2b^2\) imply
the first of the following inequalities on every finite axial
interval \(I\):
\be\l{eq:Poincare}
 \bali
 \|\bm w_h\|_{L^2(B\times I)}^2
 &\leq2\left(\|\na_h\bm w_h\|_{L^2(B\times I)}^2
                   +\|\p_3w_3\|_{L^2(B\times I)}^2\right),\\
 \|w_3\|_{L^2(B\times I)}^2
 &\leq\f4{\pi^2}\|\na_hw_3\|_{L^2(B\times I)}^2.
 \eali
\ee
The second follows from \(\int_Bw_3\d x_h=0\) and the
Poincar\'e inequality of Payne and Weinberger~\cite{PayneWeinberger}
on a convex domain of diameter two. These are also
\cite[(4.18) and (4.22)]{LiLiuZhou}. Since the derivative groups
on their right-hand sides are distinct, they imply
\be\l{eq:full-Poincare}
 \|\bm w\|_{L^2(B\times I)}^2
       \leq2\|\na\bm w\|_{L^2(B\times I)}^2.
\ee

On \(\mD_\zeta\), Young's inequality and \ef{eq:Poincare} give
\bn
 \|\bm w_h\|_{L^2(\mD_\zeta)}\|\na_hw_3\|_{L^2(\mD_\zeta)}
 &\leq\f1{2\s2}\|\bm w_h\|_{L^2(\mD_\zeta)}^2
                    +\f1{\s2}\|\na_hw_3\|_{L^2(\mD_\zeta)}^2\\
 &\leq\f1{\s2}\left(\|\na_h\bm w_h\|_{L^2(\mD_\zeta)}^2
                  +\|\p_3w_3\|_{L^2(\mD_\zeta)}^2
                  +\|\na_hw_3\|_{L^2(\mD_\zeta)}^2\right)\\
 &\leq\f1{\s2}\int_{\mD_\zeta}|\na\bm w|^2\dx.
\en
This, together with \ef{eq:centered1}, yields
\be\l{eq:basic-interaction}
 \left|\int_{\mD_\zeta}\left(P_\Phi-\f\Phi\pi\right)
                   \bm w_h\cdot\na_hw_3\dx\right|
 \leq\f{\al\Phi}{\s2(\al+4)\pi}
                   \int_{\mD_\zeta}|\na\bm w|^2\dx.
\ee

Instead of \eqref{eq:basic-interaction}, we will derive a refined estimate for \(\al>1\). By applying Proposition~\ref{prop:Robin} separately
to \(w_1\) and \(w_2\), and integrate in \(x_3\). Because
\(|\bm w_h|^2=|w_\tau|^2\) on the wall, this gives
\be\l{eq:horizontal-Robin}
 \bali
 \|\bm w_h\|_{L^2(\mD_\zeta)}^2
 \leq\f{7(\al+2)}{40(\al-1)}
 \left(\int_{\mD_\zeta}|\na_h\bm w_h|^2\dx
 +(\al-1)\int_{\p B\times(-\zeta,\zeta)}|w_\tau|^2\d S\right).
 \eali
\ee
Only horizontal integrations occur in this step; there is no
additional end-section term. Applying \ef{eq:horizontal-Robin},
Cauchy--Schwarz and then
\(ab\leq(a^2+b^2)/2\), we obtain
\bn
 &\|\bm w_h\|_{L^2(\mD_\zeta)}
                 \|\na_hw_3\|_{L^2(\mD_\zeta)}\\
 &\leq\s{\f{7(\al+2)}{40(\al-1)}}
 \left(\int_{\mD_\zeta}|\na_h\bm w_h|^2\dx
 +(\al-1)\int_{\p B\times(-\zeta,\zeta)}|w_\tau|^2\d S\right)^{1/2}
                 \|\na_hw_3\|_{L^2(\mD_\zeta)}\\
 &\leq\f12\s{\f{7(\al+2)}{40(\al-1)}}
 \left(\int_{\mD_\zeta}
        (|\na_h\bm w_h|^2+|\na_hw_3|^2)\dx
 +(\al-1)\int_{\p B\times(-\zeta,\zeta)}|w_\tau|^2\d S\right).
\en
The two bulk terms in the last line use different components of
\(\na\bm w\). Thus, by \ef{eq:centered1},
\be\l{eq:improved-interaction}
 \bali
 &\left|\int_{\mD_\zeta}\left(P_\Phi-\f\Phi\pi\right)
                  \bm w_h\cdot\na_hw_3\dx\right|\\
 &\leq\f{\al\Phi}{2(\al+4)\pi}
        \s{\f{7(\al+2)}{40(\al-1)}}
 \left(\int_{\mD_\zeta}|\na\bm w|^2\dx
 +(\al-1)\int_{\p B\times(-\zeta,\zeta)}|w_\tau|^2\d S\right).
 \eali
\ee
This is the estimate supplied by the cubic polynomial. It uses
both the bulk gradient and the circumferential boundary dissipation.

\subsection{Absorption and the proof of uniqueness}\l{Section2.3}

Throughout this subsection assume \(0\leq\Phi\leq2\pi\s{40/7}\). For \(0<\al\leq1\), the ratio between the coefficient in
\ef{eq:basic-interaction} and the dissipation coefficient in
\ef{eq:basic-dissipation} satisfies
\bes
 \left(\f{\Phi}{\s2(\al+4)\pi}\right)^2
 \leq\f{80}{7(\al+4)^2}\leq\f57<1.
\ees
For \(1\leq\al\leq3/2\), use
\(\al/(\al+4)\leq3/11\) to obtain
\bes
 \left(\f{\al\Phi}{\s2(\al+4)\pi}\right)^2
 \leq\f{80}{7}\left(\f3{11}\right)^2<1.
\ees
It follows that
\be\l{eq:small-gap}
 \min\{\al,1\}-\f{\al\Phi}{\s2(\al+4)\pi}>0\,,
 \q\text{for }0<\al\leq\f32.
\ee

For \(\al\geq3/2\), use \ef{eq:improved-interaction} and retain
the boundary term in \ef{eq:gradient}. Its coefficient satisfies
\bes
 \f{\al\Phi}{2(\al+4)\pi}\s{\f{7(\al+2)}{40(\al-1)}}
 \leq\f{\al}{\al+4}\s{\f{\al+2}{\al-1}}<1,
\ees
because
\bes
 (\al+4)^2(\al-1)-\al^2(\al+2)
 =5\al^2+8\al-16>0\,,\qq(\al\geq3/2).
\ees
In particular,
\be\l{eq:large-absorption}
 \bali
 &\int_{\mD_\zeta}\bm w\cdot\Dl\bm w\dx
 +\left|\int_{\mD_\zeta}\left(P_\Phi-\f\Phi\pi\right)
                              \bm w_h\cdot\na_hw_3\dx\right|\\
 &\leq-\left(1-\f{\al\Phi}{2(\al+4)\pi}
                     \s{\f{7(\al+2)}{40(\al-1)}}\right)
 \left(\int_{\mD_\zeta}|\na\bm w|^2\dx
 +(\al-1)\int_{\p B\times(-\zeta,\zeta)}|w_\tau|^2\d S\right)\\
 &\q-\al\int_{\p B\times(-\zeta,\zeta)}|w_3|^2\d S
 +C\int_{B\times\{\pm\zeta\}}|\bm w|\,|\na\bm w|\d x_h.
 \eali
\ee
The factor in parentheses is positive for every finite
\(\al\geq3/2\), including at the stated flux endpoint.
The remaining lateral contributions can now be discarded.
Together, \ef{eq:small-gap} and \ef{eq:large-absorption} cover
every finite \(\al>0\).

It remains to treat the Dirichlet case \(\al=\i\).
Here \(\bm w=0\) on the lateral boundary, and
\[
P_\Phi(x_h)=\f{2\Phi}{\pi}(1-|x_h|^2),
\qquad
\left\|P_\Phi-\f{\Phi}{\pi}\right\|_{L^\i(B)}
=\f{|\Phi|}{\pi}.
\]
The classical Dirichlet Poincar\'e inequality on the unit disk
gives
\[
\|g\|_{L^2(B)}
\leq \f{1}{j_{0,1}}\|\na_hg\|_{L^2(B)},
\qquad g\in H^1_0(B),
\]
where \(j_{0,1}\) is the first positive zero of the Bessel
function \(J_0\), with
\[
j_{0,1}^2\approx5.783186>\f{40}{7}.
\]
Applying this inequality to \(w_1\) and \(w_2\) on almost
every cross-section and integrating in \(x_3\), we obtain
\[
\|\bm w_h\|_{L^2(\mD_\zeta)}
\leq
\f{1}{j_{0,1}}
\|\na_h\bm w_h\|_{L^2(\mD_\zeta)}.
\]
Hence Cauchy--Schwarz and Young's inequality yield
\[
\begin{aligned}
\left|
\int_{\mD_\zeta}
\left(P_\Phi-\f{\Phi}{\pi}\right)
\bm w_h\cdot\na_hw_3\dx
\right|&\leq
\f{|\Phi|}{\pi j_{0,1}}
\|\na_h\bm w_h\|_{L^2(\mD_\zeta)}
\|\na_hw_3\|_{L^2(\mD_\zeta)}\\
&\leq
\f{|\Phi|}{2\pi j_{0,1}}
\int_{\mD_\zeta}|\na\bm w|^2\dx.
\end{aligned}
\]
By the flux assumption
\(|\Phi|\leq2\pi\sqrt{40/7}\) in
Theorem~\ref{thm:main}, we have
\[
\f{|\Phi|}{2\pi j_{0,1}}
\leq\f{\sqrt{40/7}}{j_{0,1}}<1.
\]
Set
\[
\delta=1-\f{\sqrt{40/7}}{j_{0,1}}>0.
\]
Integration by parts, using \(\bm w=0\) on the lateral
boundary, therefore gives
\[
\begin{aligned}
&\int_{\mD_\zeta}\bm w\cdot\Dl\bm w\dx
+
\left|
\int_{\mD_\zeta}
\left(P_\Phi-\f{\Phi}{\pi}\right)
\bm w_h\cdot\na_hw_3\dx
\right|\\
&\quad\leq
-\delta\int_{\mD_\zeta}|\na\bm w|^2\dx
+
\left[
\int_B\bm w\cdot\p_3\bm w\d x_h
\right]_{-\zeta}^{\zeta}.
\end{aligned}
\]
The energy identity \eqref{eq:energy} and the centering
formula \eqref{eq:centered-interaction} remain valid.
The same axial cutoff argument justifies these calculations
for weak solutions. Thus the end-slab estimates below apply
to both the Navier and Dirichlet cases, with
\(\al/(\al+4)\) interpreted as \(1\) when \(\al=\i\).

\subsubsection{End-slab estimates}

For \(K>1\), let
\bes
 \O_K^+=B\times(K-1,K),\qq
 \O_K^-=B\times(-K,-K+1).
\ees
Substitute \ef{eq:centered-interaction} into \ef{eq:energy},
apply the absorption estimates just proved, and integrate in
\(\zeta\) from \(K-1\) to \(K\). In treating the pressure term,
we first perform this integration with its signs and only then
take absolute values on the two slabs. This gives
\be\l{eq:end-slabs}
 \bali
\int_{K-1}^K\int_{\mD_\zeta}|\na\bm w|^2\dx\d\zeta\leq & C_{\al,\Phi}\Bigg[
 \int_{\O_K^+\cup\O_K^-}
   \bigl(|\bm w|\,|\na\bm w|+|\bm w|^3
                       +|P_\Phi|\,|\bm w|^2+|w_3|^2\bigr)\dx\\
 &
 +\Big|\int_{\O_K^+}(\pi-\pi_P)w_3\dx\Big|
 +\Big|\int_{\O_K^-}(\pi-\pi_P)w_3\dx\Big|\Bigg].
 \eali
\ee
The extra \(|w_3|^2\) term comes from \ef{eq:extra-end}.
In fact, its coefficient is controlled explicitly by
\bn
 \f{\al\Phi}{2(\al+4)\pi}
           \int_{\O_K^+\cup\O_K^-}|w_3|^2\dx
 &\leq\f{2\al\Phi}{(\al+4)\pi^3}
           \int_{\O_K^+\cup\O_K^-}|\na_hw_3|^2\dx\\
 &\leq\f{2\al\Phi}{(\al+4)\pi^3}
           \int_{\O_K^+\cup\O_K^-}|\na\bm w|^2\dx,
\en
where we used \ef{eq:Poincare}. Thus centering the profile
does not introduce an uncontrolled end contribution.

On either unit slab, \ef{eq:full-Poincare}, the Sobolev embedding
\(H^1\hookrightarrow L^6\), and interpolation imply
\bn
 \int_{\O_K^\pm}|\bm w|\,|\na\bm w|\dx
 &\leq C\|\na\bm w\|_{L^2(\O_K^\pm)}^2,\\
 \int_{\O_K^\pm}|\bm w|^3\dx
 &\leq\|\bm w\|_{L^2(\O_K^\pm)}^{3/2}
       \|\bm w\|_{L^6(\O_K^\pm)}^{3/2}
 \leq C\|\na\bm w\|_{L^2(\O_K^\pm)}^3,\\
 \int_{\O_K^\pm}|P_\Phi|\,|\bm w|^2\dx
 &\leq C\|P_\Phi\|_{L^\i(B)}
                  \|\na\bm w\|_{L^2(\O_K^\pm)}^2.
\en
All the slabs are translates of the same bounded domain, so the
constants are independent of \(K\).

For the pressure we use the estimate in
\cite[Proposition 4.1 and (4.15)]{LiLiuZhou}:
\be\l{eq:pressure-estimate}
 \bali
 \Big|\int_{\O_K^\pm}(\pi-\pi_P)w_3\dx\Big|
 \leq C\Big[&(1+\|P_\Phi\|_{L^\i(B)})
                    \|\na\bm w\|_{L^2(\O_K^\pm)}^2+\|\na\bm w\|_{L^2(\O_K^\pm)}^3\Big].
 \eali
\ee
This estimate is proved by solving a divergence equation on a unit
slab with right-hand side \(w_3\), whose integral is zero.
The test field has zero trace, so no lateral boundary term appears.
The argument requires no smallness of \(\Phi\) and applies to
the present perturbation without change. We therefore cite that
proof rather than repeat the same pressure calculation.

Combining \ef{eq:end-slabs} and \ef{eq:pressure-estimate} with the
preceding slab estimates yields
\be\l{eq:Saint-Venant}
 \bali
 \int_{K-1}^K\int_{\mD_\zeta}|\na\bm w|^2\dx\d\zeta
 \leq C_{\al,\Phi}\Big(&
 \|\na\bm w\|_{L^2(\O_K^+\cup\O_K^-)}^2+\|\na\bm w\|_{L^2(\O_K^+\cup\O_K^-)}^3\Big).
 \eali
\ee

\subsubsection{Completion of the proof via Ladyzhenskaya-Solonnikov approach}\l{Sec232}

	Set
	\bes
	Y(K)=\int_{K-1}^K\int_{\mD_\zeta}
	|\na\bm w|^2\dx\d\zeta,\qquad K>1.
	\ees
	Then \(Y\) is nonnegative, nondecreasing and locally
	absolutely continuous, with
	\bes
	Y'(K)=\int_{\O_K^+\cup\O_K^-}|\na\bm w|^2\dx
	\qquad\text{for a.e. }K>1.
	\ees
	Thus \ef{eq:Saint-Venant} gives
	\[
	Y(K)\leq C_{\al,\Phi}
	\bigl(Y'(K)+(Y'(K))^{3/2}\bigr)
	\qquad\text{for a.e. }K>1,
	\]
	where \(C_{\al,\Phi}\) is independent of \(K\).
	
	Since \(P_\Phi\) is independent of \(x_3\),
	the growth assumption \ef{eq:growth} implies
	\bes
	\bali
	0\leq Y(K)
	&\leq\int_{\mD_K}|\na\bm w|^2\dx\\
	&\leq2\|\na\bm u\|_{L^2(\mD_K)}^2
	+4K\|\na_hP_\Phi\|_{L^2(B)}^2
	=o(K^3).
	\eali
	\ees
	Then the following lemma, originally proved by Ladyzhenskaya and Solonnikov \cite{LadyzhenskayaSolonnikov}, implies that \(Y\equiv0\).
\begin{lemma}\l{lem:energy-growth}
Let \(Y:(1,\i)\to[0,\i)\) be nondecreasing and locally
absolutely continuous. Suppose that, for some fixed \(C>0\),
\[
 Y(K)\leq C\bigl(Y'(K)+(Y'(K))^{3/2}\bigr)
 \qquad\text{for a.e. }K>1.
\]
Then either \(Y\equiv0\) or
\bes
 \liminf_{K\to\i}\f{Y(K)}{K^3}>0.
\ees
In particular, \(Y(K)=o(K^3)\) implies \(Y\equiv0\).
\end{lemma}

	Hence \(\na\bm w=0\) in \(\mD\), and
	\ef{eq:full-Poincare} gives \(\bm w=0\).
	Thus \(\bm u=\bm P_\Phi\), as claimed.

\begin{rem}[The Theoretical Limit of the Method]
	\sl{The cubic polynomial proves a sufficient uniform flux bound larger
	than \(15\).  Within the framework of this method, this value is nearly optimal because the scalar
	Robin estimate and the absorption step used above cannot give a
	uniform sufficient flux bound of \(16\), even if we increase the degree of the polynomial in \eqref{Dpq}. Indeed, if we choose \(g(x_h)=1-|x_h|^2\) in \eqref{eq:scalar-bound}, direct calculation gives
	\bn
	\int_B|\na_hg|^2\d x_h
	&=2\pi\int_0^1 4r^3\dr=2\pi,\\
	\int_B|g|^2\d x_h
	&=2\pi\int_0^1(1-r^2)^2r\dr=\f\pi3.
	\en
and we infer from \eqref{eq:scalar-bound} that \(\la\leq6\) for every \(\al>1\).
	Using \eqref{eq:centered-interaction} and repeating the
	Cauchy--Schwarz and Young estimates leading to
	\eqref{eq:improved-interaction}, we need the following to obtain the uniqueness
\[
\Phi<
2\pi\f{\al+4}{\al}\sqrt{\lambda}
\leq
2\pi\sqrt6\left(1+\f4\al\right).
\]
Let
	\(\al\to\i\), the above upper endpoint is at most
	\bes
	2\pi\s6\approx15.3906<16.
	\ees}
\end{rem}

\section{Proof of Theorem~\ref{thm:small-friction}}\l{SecSec3}
\l{sec:proof-small}

\subsection{Energy identity and horizontal strain}

Recall that $\bm w=\bm u-\bm P_\Phi$, where
$\bm P_\Phi=P_\Phi\bm e_3$.
As in \eqref{eq:perturbation},
\begin{equation}\label{eq:general-perturbation-system}
-\Delta\bm w+P_\Phi\partial_3\bm w
+(\bm w\cdot\nabla)\bm w
+(\bm w_h\cdot\nabla_hP_\Phi)\bm e_3
+\nabla(\pi-\pi_P)=0
\quad\text{in }\mathcal D,
\end{equation}
with
\begin{equation}\label{eq:general-perturbation-constraints}
\begin{cases}
\operatorname{div}\bm w=0
&\text{in }\mathcal D,\\
\bm w_h\cdot\bm n_h=0
&\text{on }\partial\Sigma\times\mathbb R,\\
\displaystyle\int_\Sigma w_3(x_h,x_3)\,dx_h=0
&\text{for a.e. }x_3\in\mathbb R.
\end{cases}
\end{equation}
Let $\bm\tau$ be a unit tangent to $\partial\Sigma$, with
$\partial_{\bm\tau}\bm n_h=\kappa\bm\tau$, and set
$w_\tau=\bm w_h\cdot\bm\tau$.
For smooth fields, the lateral boundary conditions read
\begin{equation}\label{eq:general-component-bc}
\partial_{\bm n_h}w_\tau=(\kappa-\alpha)w_\tau,
\qquad
\partial_{\bm n_h}w_3=-\alpha w_3.
\end{equation}
Using \eqref{eq:general-component-bc}, we test
\eqref{eq:general-perturbation-system} with $\bm w$ on $\mathcal D_\zeta=\Sigma\times(-\zeta,\zeta)$,
and argue as in \eqref{eq:energy}--\eqref{eq:gradient} to obtain
\begin{equation}\label{eq:general-end-identity}
\begin{aligned}
&\int_{\mathcal D_\zeta}|\nabla\bm w|^2\,\d x
+\int_{\partial\Sigma\times(-\zeta,\zeta)}
\bigl((\alpha-\kappa)|w_\tau|^2+\alpha|w_3|^2\bigr)\,\d S
+\int_{\mathcal D_\zeta}
(\bm w_h\cdot\nabla_hP_\Phi)w_3\,\d x\\
&\quad=
\left[\int_\Sigma\left\{
\bm w\cdot\partial_3\bm w
-\frac12(P_\Phi+w_3)|\bm w|^2
-(\pi-\pi_P)w_3
\right\}\,\d x_h\right]_{-\zeta}^{\zeta}.
\end{aligned}
\end{equation}
For weak solutions, \eqref{eq:general-end-identity} holds
for almost every $\zeta>0$ by the same axial cutoff argument
used for \eqref{eq:energy}. 

For $\bm{v}\in H^1(\Sigma,\mR^2)$, we denote
\[
\mathbb S_h\bm v
\ed\frac12\bigl(\nabla_h\bm v+
(\nabla_h\bm v)^{\mathsf T}\bigr).
\]
Thus for smooth tangential fields, expanding the symmetric gradient gives
\be\l{Int1021-0}
2\int_\Sigma |\mathbb{S}_h\boldsymbol{w}_h|^2\,\dx_h
=
\int_\Sigma |\nabla_h\boldsymbol{w}_h|^2\,\dx_h
+
\sum_{i,j=1}^2
\int_\Sigma \partial_i w_j\,\partial_j w_i\,\dx_h.
\ee
Integrating the last term by parts first in $x_i$, we obtain
\be\l{Int1021}
\begin{aligned}
\sum_{i,j=1}^2
\int_\Sigma \partial_i w_j\,\partial_j w_i\,\dx_h
&=
\sum_{i,j=1}^2
\int_{\partial\Sigma} n_i w_j\,\partial_j w_i\,\d S_h
-
\sum_{i,j=1}^2
\int_\Sigma w_j\,\partial_i\partial_j w_i\,\dx_h
\\
&=
\un{\int_{\partial\Sigma}
\bigl((\boldsymbol{w}_h\cdot\nabla_h)\boldsymbol{w}_h\bigr)
\cdot\boldsymbol{n}_h\,\d S_h}_{B_1}
-
\int_\Sigma
\boldsymbol{w}_h\cdot\nabla_h
(\operatorname{div}_h\boldsymbol{w}_h)\,\dx_h.
\end{aligned}
\ee
A second integration by parts gives
\be\l{Int1021-1}
\begin{aligned}
-\int_\Sigma
\boldsymbol{w}_h\cdot\nabla_h
(\operatorname{div}_h\boldsymbol{w}_h)\,\dx_h
&=
\int_\Sigma
|\operatorname{div}_h\boldsymbol{w}_h|^2\,\dx_h-
\int_{\partial\Sigma}
(\operatorname{div}_h\boldsymbol{w}_h)
(\boldsymbol{w}_h\cdot\boldsymbol{n}_h)\,\d S_h
\\
&=
\int_\Sigma
|\operatorname{div}_h\boldsymbol{w}_h|^2\,\dx_h,
\end{aligned}
\ee
where the boundary term vanishes because
$\boldsymbol{w}_h\cdot\boldsymbol{n}_h=0$ on $\partial\Sigma$. To compute the remaining boundary term $B_1$ in \eqref{Int1021}, write
$\boldsymbol{w}_h=w_\tau\boldsymbol{\tau}$ on $\partial\Sigma$
and recall that
$\partial_\tau\boldsymbol{n}_h=\kappa\boldsymbol{\tau}$.
Differentiating the tangency condition along the boundary yields
\be\l{Int1021-2}
\begin{aligned}
0
&=
w_\tau\partial_\tau
(\boldsymbol{w}_h\cdot\boldsymbol{n}_h)=
\bigl((\boldsymbol{w}_h\cdot\nabla_h)\boldsymbol{w}_h\bigr)
\cdot\boldsymbol{n}_h
+
w_\tau\boldsymbol{w}_h\cdot
\partial_\tau\boldsymbol{n}_h
\\
&=
\bigl((\boldsymbol{w}_h\cdot\nabla_h)\boldsymbol{w}_h\bigr)
\cdot\boldsymbol{n}_h
+
\kappa|w_\tau|^2.
\end{aligned}
\ee
Substituting \eqref{Int1021-1}--\eqref{Int1021-2} in \eqref{Int1021}, then one infers from \eqref{Int1021-0} that
\begin{equation}\label{eq:general-horizontal-identity}
2\int_\Sigma |\mathbb{S}_h\boldsymbol{w}_h|^2\,\dx_h
=
\int_\Sigma
\left(
|\nabla_h\boldsymbol{w}_h|^2
+
|\operatorname{div}_h\boldsymbol{w}_h|^2
\right)\,\dx_h
-
\int_{\partial\Sigma}\kappa|w_\tau|^2\,\d S_h.
\end{equation}
By density, \eqref{eq:general-horizontal-identity} also holds
for tangential $H^1(\Sigma;\mathbb{R}^2)$ fields.

Set
\be\l{DEFEZ}
E_\zeta
\ed\int_{\mathcal D_\zeta}|\nabla\bm w|^2\,dx
+\int_{\partial\Sigma\times(-\zeta,\zeta)}
\bigl((\alpha-\kappa)|w_\tau|^2+\alpha|w_3|^2\bigr)\,\d S.
\ee
Using the divergence-free property of $\bm{w}$, one has $\partial_3w_3=-\operatorname{div}_h\bm w_h$. Thus
\[
|\nabla\bm w|^2
=|\nabla_h\bm w_h|^2+|\nabla_hw_3|^2
+|\partial_3\bm w_h|^2
+|\operatorname{div}_h\bm w_h|^2.
\]
Substituting this in \eqref{DEFEZ} and using \eqref{eq:general-horizontal-identity}, the two curvature terms cancel. Since
$\bm w_h\cdot\bm n_h=0$ on the lateral boundary,
$|\bm w|^2=|w_\tau|^2+|w_3|^2$ there. Hence
\begin{equation}\label{eq:general-dissipation}
	\begin{aligned}
		E_\zeta
		={}&\int_{\mathcal D_\zeta}
		\bigl(2|\mathbb S_h\bm w_h|^2
		+|\nabla_hw_3|^2+|\partial_3\bm w_h|^2\bigr)\,\d x+\alpha\int_{\partial\Sigma\times(-\zeta,\zeta)}
		|\bm w|^2\,\d S
		\geq0.
	\end{aligned}
	\end{equation}
This identity incorporates the curvature term into the horizontal
strain. Its bulk terms remain unweighted by $\alpha$, which is
essential for absorbing the interaction at small friction.

Set
\be\l{DIzeta}
I_\zeta\ed\int_{\mathcal D_\zeta}
(\bm w_h\cdot\nabla_hP_\Phi)w_3\,\d x.
\ee
As in the previous section, we will use $E_\zeta$ to eliminate $I_\zeta$. We will divide the detailed argument into two cases: the non-disk cross-section case and the disk cross-section case.

\subsection{Non-disk cross-section case}

We first introduce the following form of Korn's inequality.

\begin{lemma}[Korn inequality on non-disk cross-sections]
\label{lem:noncircular-Korn}
Let $\Sigma\subset\mathbb R^2$ be a bounded, smooth, simply
connected domain that is not a disk. Then
\begin{equation}\label{eq:noncircular-Korn}
\|\bm v\|_{H^1(\Sigma)}
\leq C_\Sigma\|\mathbb S_h\bm v\|_{L^2(\Sigma)}
\end{equation}
for every $\bm v\in H^1(\Sigma;\mathbb R^2)$ satisfying
$\bm v\cdot\bm n_h=0$ on $\partial\Sigma$,
where $C_\Sigma$ depends only on $\Sigma$.
\end{lemma}

\begin{proof}
If $\mathbb S_h\bm v=0$, then
{$\bm v(x_h)=\bm a+bx_h^\perp$} for some constant vector
$\bm a$ and $b\in\mathbb R$. We \textbf{claim} that the boundary condition and the assumption that
$\Sigma$ is not a disk force $\bm a=0$ and $b=0$.

Suppose now that $b\neq0$. Set
\[
x_0=\left(-\frac{a_2}{b},\,\frac{a_1}{b}\right).
\]
Then
\[
\bm v(x_h)
=\bm a+b(-x_2,x_1)
=b\bigl(-(x_2-x_{0,2}),\,x_1-x_{0,1}\bigr),
\]
so $\bm v$ is a rigid rotation about $x_0$.

On $\partial\Sigma$, write $\bm n_h=(n_1,n_2)$ and choose
the unit tangent $\bm\tau=(-n_2,n_1)$. The boundary condition gives
\[
0=\bm v\cdot\bm n_h
=-b(x_h-x_0)\cdot\bm\tau.
\]
Since $b\neq0$, it follows that
\[
\partial_\tau|x_h-x_0|^2
=2(x_h-x_0)\cdot\bm\tau=0\,,
\qquad\text{on }\partial\Sigma.
\]
Thus the distance from $x_0$ is constant along each connected
component of $\partial\Sigma$. Since $\Sigma$ is bounded, smooth and simply connected, $\partial\Sigma$ is a smooth closed curve, so it must be the entire circle. Thus
\[
\Sigma=\{x_h\in\mathbb R^2:|x_h-x_0|<R\},
\]
contrary to the assumption that $\Sigma$ is not a disk.
Therefore $b=0$. The boundary condition $\bm v\cdot\bm n_h=0$ then implies $\bm{a}=0$. This proves the claim.

Now we prove \eqref{eq:noncircular-Korn}. Suppose it fails, then there exists
a sequence $\bm v_j\in H^1(\Sigma;\mathbb R^2)$ such that
\[
\bm v_j\cdot\bm n_h=0
\quad\text{on }\partial\Sigma,
\qquad
\|\bm v_j\|_{H^1(\Sigma)}=1,
\qquad
\|\mathbb S_h\bm v_j\|_{L^2(\Sigma)}\rightarrow0.
\]
Since $\{\bm v_j\}$ is bounded in $H^1(\Sigma;\mathbb R^2)$,
weak compactness and Rellich compactness yield, after passing
to a subsequence,
\[
\bm v_j\rightharpoonup\bm v
\quad\text{in }H^1(\Sigma;\mathbb R^2),
\qquad
\bm v_j\rightarrow\bm v
\quad\text{in }L^2(\Sigma;\mathbb R^2).
\]
By continuity of the trace operator, the limit satisfies
\[
\bm v\cdot\bm n_h=0
\qquad\text{on }\partial\Sigma.
\]
Moreover, weak convergence in $H^1$ implies
\[
\mathbb S_h\bm v_j\rightharpoonup\mathbb S_h\bm v
\quad\text{in }L^2(\Sigma).
\]
Since $\|\mathbb S_h\bm v_j\|_{L^2(\Sigma)}\to0$, by applying the weak lower semicontinuity of $L^2$ norm, we obtain $\mathbb S_h\bm v=0$. Thus $\bm v=\bm a+b(-x_2,x_1)$, and the preceding claim
gives $\bm v=0$, so the strong
$L^2$ convergence yields
\[
\|\bm v_j\|_{L^2(\Sigma)}\longrightarrow0.
\]
Korn's second inequality (c.f. \cite[Theorem 2]{BauerPauly}) then gives
\[
1=\|\bm v_j\|_{H^1(\Sigma)}
\leq C_\Sigma
\left(
\|\mathbb S_h\bm v_j\|_{L^2(\Sigma)}
+\|\bm v_j\|_{L^2(\Sigma)}
\right)
\longrightarrow0,
\]
which is a contradiction. This completes the proof of the lemma.
\end{proof}

By \cite[Proposition~3.1]{LiLiuZhou},
\begin{equation}\label{eq:general-profile-gradient}
\|\nabla_hP_\Phi\|_{L^2(\Sigma)}
\leq C_\Sigma\frac{\alpha|\Phi|}{1+\alpha},
\qquad\alpha>0.
\end{equation}
Lemma~\ref{lem:noncircular-Korn}, Sobolev embedding and
the zero mean of $w_3$ give
\[
\|\bm w_h\|_{L^4(\Sigma)}
\leq C_\Sigma\|\mathbb S_h\bm w_h\|_{L^2(\Sigma)},
\qquad
\|w_3\|_{L^4(\Sigma)}
\leq C_\Sigma\|\nabla_hw_3\|_{L^2(\Sigma)}.
\]
Recalling the definition \eqref{DIzeta}, H\"older's and Young's inequalities, together with
\eqref{eq:general-profile-gradient}, yield
\begin{equation}\label{eq:general-direct-interaction}
\begin{aligned}
|I_\zeta|
&\leq C_\Sigma\frac{\alpha|\Phi|}{1+\alpha}
\int_{-\zeta}^{\zeta}
\|\mathbb S_h\bm w_h\|_{L^2(\Sigma)}
\|\nabla_hw_3\|_{L^2(\Sigma)}\,\dx_3\\
&\leq C_\Sigma\frac{\alpha|\Phi|}{1+\alpha}
\int_{\mathcal D_\zeta}
\bigl(2|\mathbb S_h\bm w_h|^2+|\nabla_hw_3|^2\bigr)\,\dx\\
&\leq C_\Sigma\frac{\alpha|\Phi|}{1+\alpha}E_\zeta.
\end{aligned}
\end{equation}
Since $\partial_3 w_3=-\operatorname{div}_h\bm w_h$,
\eqref{eq:general-dissipation} and \eqref{eq:noncircular-Korn} give
\[
\|\nabla\bm w\|_{L^2(\mathcal D_\zeta)}^2
\leq C_\Sigma E_\zeta.
\]
Choosing $c_\Sigma$ in Theorem~\ref{thm:small-friction} properly
so that $C_\Sigma c_\Sigma\leq1/2$, we obtain from
\eqref{eq:general-direct-interaction}
\begin{equation}\label{eq:general-absorption}
E_\zeta-|I_\zeta|
\geq\frac12E_\zeta
\geq C_\Sigma^{-1}
\int_{\mathcal D_\zeta}|\nabla\bm w|^2\,\dx,
\end{equation}
after enlarging $C_\Sigma$ if necessary.

\subsection{Disk cross-section case}

Suppose that $\Sigma=B$. Since the profile
\eqref{eq:Poiseuille} is radial, we first estimate
$x_h\cdot\bm w_h$.

Set $r\ed|x_h|$, $\bm{e_r}\ed\frac{x_h}{|x_h|}= (\cos\theta,\sin\theta)$ and $w_r=\bm w_h\cdot\bm e_r$.
By \eqref{eq:general-perturbation-constraints},
$w_r(1)=0$. For smooth fields,
\[
\partial_rw_r
=\bm e_r^{\mathsf T}(\nabla_h\bm w_h)\bm e_r
=\bm e_r^{\mathsf T}\mathbb S_h\bm w_h\bm e_r.
\]
For each fixed $(\theta,x_3)$, the Cauchy--Schwarz
inequality gives
\[
\begin{aligned}
	|w_r(r)|^2
	=\left|\int_r^1\partial_sw_r(s)\,\d s\right|^2&\leq
	\left(\int_r^1s|\partial_sw_r(s)|^2\,\d s\right)
	\left(\int_r^1\frac{\d s}{s}\right)\\
	      &\leq
	\log(1/r)\int_0^1s|\partial_sw_r(s)|^2\,\d s.
\end{aligned}
\]
Multiplying by $r^3$ and integrating, we obtain
\[
\int_0^1r^3|w_r(r)|^2\,\d r
\leq\frac1{16}\int_0^1r|\partial_rw_r(r)|^2\,\d r,
\]
where $\int_0^1r^3\log(1/r)\,\d r=1/16$.
Since $x_h\cdot\bm w_h=rw_r$ and
$|\partial_rw_r|\leq|\mathbb S_h\bm w_h|$,
integration in $\theta$ and $x_3$ yields
\begin{equation}\label{eq:small-radial-estimate}
	\|x_h\cdot\bm w_h\|_{L^2(\mathcal D_\zeta)}
	\leq\frac14
	\|\mathbb S_h\bm w_h\|_{L^2(\mathcal D_\zeta)}.
\end{equation}
For $H^1$ fields, the same argument applies on annuli,
followed by letting the inner radius tend to zero and the density argument.

Using \eqref{eq:Poiseuille},
\eqref{eq:small-radial-estimate} and
\eqref{eq:Poincare}, followed by Young's inequality, we obtain
\begin{equation}\label{eq:small-direct-interaction1}
	\begin{aligned}
		|I_\zeta|
		&=\frac{4\alpha|\Phi|}{\pi(\alpha+4)}
		\left|\int_{\mathcal D_\zeta}
		(x_h\cdot\bm w_h)w_3\,\d x\right|\\
		&\leq\frac{2\alpha|\Phi|}{\pi^2(\alpha+4)}
		\|\mathbb S_h\bm w_h\|_{L^2(\mathcal D_\zeta)}
		\|\nabla_hw_3\|_{L^2(\mathcal D_\zeta)}\\
		&\leq\frac{\alpha|\Phi|}{\sqrt2\,\pi^2(\alpha+4)}
		\int_{\mathcal D_\zeta}
		\bigl(2|\mathbb S_h\bm w_h|^2
		+|\nabla_hw_3|^2\bigr)\,\d x.
	\end{aligned}
\end{equation}
Noticing that
\[
\frac{\alpha|\Phi|}
{\sqrt2\,\pi^2(\alpha+4)}<1,
\]
which follows from
\eqref{eq:small-friction-threshold}, by \eqref{eq:general-dissipation} and
\eqref{eq:small-direct-interaction1}, using
$|\bm w|^2\geq|w_\tau|^2$ on the lateral boundary, we obtain
\begin{equation}\label{eq:small-coercivity}
\begin{aligned}
E_\zeta-|I_\zeta|
&\geq
(1-\frac{\alpha|\Phi|}
{\sqrt2\,\pi^2(\alpha+4)})\int_{\mathcal D_\zeta}
\bigl(
2|\mathbb S_h\bm w_h|^2
+|\nabla_hw_3|^2
+|\partial_3\bm w_h|^2
\bigr)\,\d x\\
&\q+
\alpha\int_{\partial B\times(-\zeta,\zeta)}
|w_\tau|^2\,\d S
\\
&\geq
\min\left\{
\alpha,\,
1-\frac{\alpha|\Phi|}{\sqrt2\,\pi^2(\alpha+4)}
\right\}
\int_{\mathcal D_\zeta}|\nabla\bm w|^2\,\d x.
\end{aligned}
\end{equation}
Here the last step holds owing to
\[
\int_{\mathcal D_\zeta}|\nabla\bm w|^2\,\d x
=
\int_{\mathcal D_\zeta}
\bigl(
2|\mathbb S_h\bm w_h|^2
+|\nabla_hw_3|^2
+|\partial_3\bm w_h|^2
\bigr)\,\d x
+
\int_{\partial B\times(-\zeta,\zeta)}|w_\tau|^2\,\d S\,,
\]
which is a direct conclusion by combining \eqref{DEFEZ} and \eqref{eq:general-dissipation}, and noticing $\kappa=1$ on $\partial B$.

\begin{rem}[Cylindrical pipes of arbitrary radius]
	Let $\Sigma=B_R(x_0)$ be a disk of radius $R>0$, and let
	$\alpha>0$. If $\bm u$ is a weak solution of
	\eqref{eq:NS}--\eqref{eq:flux} satisfying \eqref{eq:growth}, then
	$\bm u=\bm P_\Phi$ whenever
	\[
	|\Phi|<
	\frac{\sqrt2\,\pi^2(R\alpha+4)}{\alpha}.
	\]
	In particular, for circular cross-sections, one may choose $c_\Sigma=\sqrt2\,\pi^2\min\{1,R\}$ in Theorem~\ref{thm:small-friction}.
\end{rem}

\subsection{Completion of the proof}\label{sec:small-completion}

In both cases, \eqref{eq:general-absorption} and
\eqref{eq:small-coercivity} give
\[
E_\zeta-|I_\zeta|
\geq c_{\Sigma,\alpha,\Phi}
\int_{\mathcal D_\zeta}|\nabla\bm w|^2\,\dx,
\qquad c_{\Sigma,\alpha,\Phi}>0,
\]
with a constant independent of $\zeta$.
We now argue in Section \ref{Section2.3}, with $B$ replaced by $\Sigma$,
to obtain $\bm w=0$. Hence $\bm u=\bm P_\Phi$.

\end{document}